\documentclass[centertags,12pt]{amsart}
\usepackage{latexsym}
\usepackage{amssymb}

\numberwithin{equation}{section}
\makeatletter
\renewcommand{\subsection}{\@startsection
{subsection}{2}{0mm}{\baselineskip}{-0.25cm}
{\normalfont\normalsize\bf}}
\makeatother

\newtheorem{theorem}{Theorem}[section]

\newtheorem{lemma}[theorem]{Lemma}

{\theoremstyle{definition}
\newtheorem{example}[theorem]{Example}}
\theoremstyle{remark}
\newtheorem{remark}[theorem]{Remark}

\def\P{\mathbf P}
\def\Z{\mathbf Z}
\def\cC{\mathcal C}
\def\cD{\mathcal D}
\def\cH{\mathcal H}
\def\cJ{\mathcal J}
\def\cL{\mathcal L}
\def\cR{\mathcal R}
\def\cX{\mathcal X}
\def\cY{\mathcal Y}
\def\cZ{\mathcal Z}
\def\fl{\mathbf{F_\ell}}
\def\fq{\mathbf F_{q^2}}
\def\dim{{\rm dim}}
\def\div{{\rm div}}
\def\deg{{\rm deg}}
\def\det{{\rm det}}

\def\frx{{\mathbf{Fr}}_{\mathcal X}}
\def\frj{{\mathbf{Fr}}_{\mathcal J}}
\def\gep{\epsilon}

\begin{document}

\author[G.~Korchm\'aros]{G\'abor Korchm\'aros}
\author[F.~Torres]{Fernando Torres}\thanks{This research was carried out
within the activity of GNSAGA of the Italian CNR with the support of the
Italian Ministry for Research and Technology. The present paper was partly
written while the second author was visiting ICTP-Italy (July-August 1999),
as a Regular Associate, and the University of Basilicata at Potenza-Italy
(July 1999)}

\title[Maximal curves in a hermitian variety]{Embedding of a maximal curve
in a\\
hermitian variety}

   \address{Dipartimento di Matematica, Universit\'a della
Basilicata, via N. Sauro 85,\\
85100 Potenza, Italy}
\email{korchmaros@unibas.it}
   \address{IMECC-UNICAMP, Cx. P. 6065, Campinas-13083-970-SP, Brazil}
\email{ftorres@ime.unicamp.br}

\begin{abstract}
                Let $\cX$ be a projective geometrically irreducible
non-singular algebraic curve defined over a finite field $\fq$ of order
$q^2$. If the number of $\fq$-rational points of $\cX$ satisfies the
Hasse-Weil upper bound, then $\cX$ is said to be $\fq$-maximal. For a point
$P_0\in \cX(\fq)$, let $\pi$ be the morphism arising from the linear series
$\cD:=|(q+1)P_0|$, and let $N:=\dim(\cD)$. It is known that $N\ge 2$ and
that $\pi$ is independent of $P_0$ whenever $\cX$ is $\fq$-maximal. The
following theorems will be proved:

\begin{theorem}\label{thmA} If $\cX$ is $\fq$-maximal, then $\pi:\cX\to
\pi(\cX)$ is a $\fq$-isomorphism. The non-singular model $\pi(\cX)$ has
degree $q+1$ and lies on a Hermitian variety defined over $\fq$ of
$\P^N(\bar\fq)$. 
\end{theorem}

\begin{theorem}\label{thmB} If $\cX$ is $\fq$-maximal, then it is
$\fq$-isomorphic to a curve $\cY$ in $\P^M(\bar\fq)$, with $2\le M\le N$,
such that $\cY$ has degree $q+1$ and lies on a non-degenerate Hermitian
variety defined over $\fq$ of $\P^M(\bar\fq)$. Furthermore, ${\rm
Aut}_{\fq}(X)$ is isomorphic to a subgroup of the projective unitary group
$PGU(M+1,q^2)$. 
\end{theorem}

\begin{theorem}\label{thmC} If $\cX$ is $\fq$-birational to a curve $\cY$
embedded in $\P^M(\bar\fq)$ such that $\cY$ has degree $q+1$ and lies on a
non-degenerate Hermitian variety defined over $\fq$ of $\P^M(\bar\fq)$,
then $\cX$ is $\fq$-maximal and $\cX$ is $\fq$-isomorphic to $\cY$. 
\end{theorem}
    \end{abstract}

\maketitle

\section{Introduction}\label{s1}

Let $\cX$ be a projective geometrically irreducible non-singular algebraic
curve defined over $\fl$, the finite field of order $\ell$. There is a
natural way to define a $\fl$-linear series $\cD$ on the curve $\cX$
provided that $\cX(\fl)\neq\emptyset$, and geometrical and arithmetical
properties of $\cX$ may be investigated by using $\cD$.  This linear
series $\cD$ arises from the characteristic polynomial $h(t)$ of the
Jacobian $\cJ$ (over $\fl$) of $\cX$ in the following way, see
\cite[Section 1.3]{FT1}.  Let $\prod_{i=1}^{T}h^{r_i}_i(t)$ be the
factorization of $h(t)$ over $\Z[t]$.  Since the Frobenius morphism $\frj$
(over $\fl$) on $\cJ$ is semisimple and the representation of
endomorphisms of $\cJ$ on the Tate module is faithful \cite[Theorem
2]{Tate}, \cite[VI, Section 3]{Lang}, we have
   \begin{equation*}
\prod_{i=1}^{T} h_i(\frj)=0\qquad \text{on\quad $\cJ$}\, .
   \end{equation*}
Now let $P_0\in \cX(\fl)$ and set $m:=|\prod_{i=1}^{T}h_i(1)|$. Then the
foregoing equation is equivalent to the following linear equivalence of
$\fl$-divisors on $\cX$:
   \begin{equation}\label{eq1.1}
\sum_{i=1}^{T}\alpha_i\frx^{T-i}(P)+\frx(P)\sim mP_0\qquad P\in \cX\, ,
    \end{equation}
where $\sum_{i=1}^{T}\alpha_it^{T-i}+t^T:=\prod_{i=1}^{T}h_i(t)$; see
\cite[Section 1.3]{FT1}.

Assume from now on that $\ell$ is a square, and let $q:=\sqrt \ell$. Then
$h(t)=(t+q)^{2g}$ if and only if $\cX$ is $\fq$-maximal, that is
$\#\cX(\fq)$ attains the Hasse-Weil upper bound
   $$
1+q^2+2qg\, ,
   $$
where $g$ is the genus of $\cX$. From (\ref{eq1.1}), every
$\fq$-maximal curve $\cX$ is equipped with a $\fq$-linear series
$\cD=\cD_\cX=|(q+1)P_0|$ which is independent of $P_0\in
\cX(\fq)$ and satisfies the so-called ``Fundamental Equivalence"
\cite [Corollary 1.2]{FGT}:
      \begin{equation}\label{eq1.2}
qP+\frx(P)\sim (q+1)P_0\qquad \text{for any $P\in\cX$}\, .
      \end{equation}
In particular, $(q+1)P\sim (q+1)P_0$ for all points $P\in\cX(\fq)$, see
\cite[Lemma 1]{RS}.

Maximal curves have been intensively studied also in connection with
coding theory and cryptography. The pioneer work by St\"ohr and Voloch
\cite{SV}, giving among other things an alternative proof of the
Hasse-Weil bound via Weierstrass Point Theory, has been widely used to
investigate maximal curves, their $\cD$-Weierstrass points and the support
of the $\fq$-Frobenius divisor associated to $\cD$. However, the
fundamental question in this context, namely whether the $\fq$-morphism
$\pi:\cX\to\pi(\cX)$ associated to $\cD$ is an isomorphism, has only had a
partial answer so far \cite[Proposition 1.9]{FGT}. Our Theorem
\ref{thm2.1} (which is the first statement in Theorem \ref{thmA}) states
that $\pi$ is indeed an isomorphism. (This result was originally stated in
\cite[Section 2.3]{FT1} but the proof giving there was not correct.) Hence
the maximal curve $\cX$ may be identified with a curve of degree $q+1$
embedded in the projective space $\P^N(\bar\fq)$ with $N=\dim(\cD)$. 

This allows to investigate in some more details the geometric behaviour of
$\cX$. In the smallest case, $N=2$, $\cX$ is a non-degenerate Hermitian
curve, according to the R\"uck-Stichtenoth theorem, see \cite{RS}. Our
Theorem \ref{thm3.1} (which is actually the second statement in Theorem
\ref{thmA}) is a generalization for $N>2$ as it states that $\cX$ lies on
a Hermitian variety $\cH\subseteq \P^N(\bar\fq)$ defined over $\fq$. It
might be that $\cH$ is degenerate in some cases, such a possibility
occuring when $\cX$ is $(N-1)$-strange that is the osculating hyperplanes
to $\cX$ at generic points have a non-empty intersection. This kind of
pathology in positive characteristic has been considered by several
authors after Lluis and Samuel, see the most recent papers \cite{K},
\cite{HM}, \cite{GV3}, \cite{HK} on this subject. What we are able to
prove in this direction is the existence of a projection $\phi:
\P^N(\bar\fq) \rightarrow \P^M(\bar\fq)$ such that $\phi(\pi(\cX))$ lies
on a {\em non-degenerate} Hermitian variety defined over $\fq$ of
$\P^M(\bar\fq)$, see Theorem \ref{thmB} and Section \ref{s3}. Here $M$ is
the dimension of the smallest linear series $\cR$ containing all divisors
$qP+\frx(P)$ with $P$ ranging over $\cX$. In other words, $M$ is the
dimension of the smallest $\fq$-vector subspace $V$ of the function field
$\fq(\cX)$ such that for any two points $P_1,P_2\in\cX$ there exists $f\in
V$ satisfying $qP_1+\frx(P_1)=qP_2+\frx(P_2)+\div(f)$. The converse of the
first statement of Theorem \ref{thmB} also holds, see Theorem \ref{thmC}
and Section \ref{s4}. Putting together these two theorems we see that the
study of $\fq$-maximal curves is equivalent to that of projective
geometrically irreducible non-singular curves of degree $q+1$ lying on a
non-degenerate Hermitian variety defined over $\fq$ in a projective space
over $\bar\fq$. Note that $q+1$ is the minimum degree that a non-singular
curve of degree bigger than one lying on a non-degenerate Hermitian
variety can have.

\section{Maximal curves and their natural embedding in a projective
space}\label{s2}

Our terminology in this and subsequent sections is the same as employed
in section $2$ of \cite{SV}, and in \cite{FGT}.

In this section we assume that $\cX$ is a $\fq$-maximal curve. Our aim is
to show that $\cX$ is $\fq$-isomorphic to a curve in $\P^N(\bar\fq)$,
$N$ being the dimension of the linear series $\cD=\cD_\cX=|(q+1)P_0|$ with
$P_0\in\cX(\fq)$. Let $\pi: \cX \rightarrow \P^N(\bar\fq)$ be the
morphism associated to $\cD$. Then

\begin{lemma}\label{lemma2.1} {\rm (\cite[Prop. 1.9]{FGT})} The
following statements are equivalent:
\begin{enumerate}
\item[\rm(1)] $\cX$ is $\fq$-isomorphic to $\pi(\cX);$
\item[\rm(2)] For $P\in \cX$, $\pi(P)\in
\P^N(\fq)\Leftrightarrow P\in \cX(\fq);$
\item[\rm(3)] For $P\in \cX$, $q$ is a Weierstrass non-gap at $P.$
\end{enumerate}
    \end{lemma}

Hence we can limit ourselves to prove the above statement (2). To do this we
need some previous results concerning $\cD$-orders and Frobenius $\cD$-orders.

Let $\epsilon_0=0<\epsilon_1=1<\ldots<\epsilon_N$ (resp.
$\nu_0=0<\nu_1<\ldots<\nu_{N-1}$) denote the $\cD$-orders (resp. the
$\fq$-Frobenius orders) of $\cD$.

    \begin{lemma}\label{lemma2.2} {\rm (\cite[Thm. 1.4]{FGT})} The
$\cD$-orders and Frobenius $\cD$-orders of $\cX$ have the following
properties:

    \begin{enumerate}
\item[\rm(1)] $\epsilon_N=q;$
\item[\rm(2)] $\nu_{N-1}=q;$
\item[\rm(3)] $\nu_1=1$ iff $N\ge 3;$
\item[\rm(4)] $0,1$, and $q$ (resp. $q+1$) are $(\cD,P)$-orders provided
that $P\not\in\cX(\fq)$ (resp. $P\in \cX(\fq)$)$.$
     \end{enumerate}
    \end{lemma}

Let  $\pi = (f_0:\ldots:f_N)$ where each projective coordinate $f_i$ belongs
to $\fq(\cX)$, the function field over $\fq$ of $\cX$. As in \cite {SV},
we will consider $\pi: \cX \rightarrow \P^N{(\bar\fq)}$ as a parametrized
curve in $\P^N({\bar\fq})$, and the points $P\in\cX$ will be viewed as its
places (or branches). Then the intersection divisor $\pi^{-1}(H)$
of $\cX$ arising from a hyperplane $H$ of $\P^N({\bar\fq})$ is defined in the usual
manner, and $\cD$ turns out be the linear series of hyperplane sections,
see \cite[p.3]{SV}. In particular, the osculating hyperplane at
$P$ is the hyperplane in $\P^N({\bar\fq})$ which intersects the
branch $P$ with multiplicity $j_N$, where $(j_0,j_1,\ldots,j_N)$ is
the $(\cD,P)$-order sequence, see \cite[p.4]{SV}.

Put $\cL((q+1)P_0)=\langle f_0,\ldots,f_N\rangle$. By Lemma
\ref{lemma2.2}(1) and \cite[Thm.1]{GV}, there exist $z_0,\ldots,z_N\in
\fq(\cX)$, not all zero, such that

   \begin{equation}\label{eq2.1}
z_0^qf_0+\ldots+z_N^qf_N=0\, .
   \end{equation}

Some features of the homogeneous $N$-tuple $(z_0,\ldots,z_N)$ are
stated in the following lemma.

   \begin{lemma}\label{lemma2.3}
   \begin{enumerate}
   \item[\rm(1)] The osculating hyperplane at $P\in \cX$ has equation
$$
w_0^q(P) X_0+w_1^q(P)X_1+\ldots+w_N^q(P)X_N = 0\, ,
$$
where $w_i:=t^{e_P}z_i$, $t$ a local parameter at $P$, and
$e_P:=-\min\{v_P(z_0),\ldots,v_P(z_N)\};$

    \item[\rm(2)] The following relation also holds:

\begin{equation}\label{eq2.2}
z_0f_0^q+\ldots+z_Nf^q_N=0\, ;
\end{equation}

   \item[\rm(3)] The $\fq$-rational functions $z_0,z_1,\ldots z_N$ are
uniquely determined by Eq.~(\ref{eq2.1}) up to a non-zero factor in
$\fq(\cX);$
   \end{enumerate}
   \end{lemma}
    \begin{proof} (1) For $i=0,\ldots,N$, let
$$
w_i(t)= \sum_{j=0}^{\infty}a^{(i)}_jt^j\in \bar\fq[[t]]
$$
be the local expansion of $w_i$ at $P$. As there exists $i\in
\{0,\ldots,N\}$ such that $a^{(i)}_0\neq 0$ (e.g. $i$ satisfying
$e_P=-v_P(z_i)$), we can consider the following hyperplane in
$\P^N(\bar\fq)$:
   $$
H:\qquad \sum (a^{(i)}_0)^qX_i = 0\, .
   $$
Then, thanks to Lemma \ref{lemma2.2}(4), Item (1) follows, once we have
shown that $v_P(\pi^{-1}(H))\ge q$. Taking into consideration
Eq.~(\ref{eq2.1}),
  \begin{equation}\label{eq2.3}
v_P(\sum_{i=0}^{N}
(a^{(i)}_0)^qf_i)=v_P(t^q\sum_{i=0}^{N}\sum_{j=1}^{\infty}
a^{(i)}_jt^{qj-q}f_i)\, ,
  \end{equation}
yielding the desired relation $v_P(\pi^{-1}(H))\ge q$.

(2) By the Fundamental Equivalence (\ref{eq1.2}), $\frx(P)$ belongs to the
osculating hyperplane at $P$ for every $P\in \cX$. Then from Eq.
(\ref{eq2.1}) we infer for all but a finitely many points $P\in \cX$ that
$$
\sum_{i=0}^{N}z_i(P)^qf_i(P)^{q^2}=0
$$
and Item (2) follows.

(3) This is clear because once the projective coordinates are fixed, then
the osculating hyperplane at any point is uniquely determined modulo a
non-zero element of $\bar\fq$.
    \end{proof}

        \begin{lemma}\label{lemma2.4} Let
$P\in\cX$ be such that $\pi(P)\in \P^N(\fq)$. Then $P\in\cX(\fq).$
        \end{lemma}

        \begin{proof} Since $\pi(P)$ is $\fq$-rational we can take it to
the point $(1:0:\ldots:0)$ by means of a $\fq$-linear transformation. The
new coordinates still satisfy Eqs. (\ref{eq2.1}) and (\ref{eq2.2}). In
addition, we can assume that $\pi=(1:f_1:\ldots:f_N)$ so that $v_P(f_i)\ge
1$ for $i\ge 1$. Now, the set up and the results of the computation
involving local expansion in the proof of Lemma \ref{lemma2.3}(2) together
with Lemma \ref{lemma2.2}(4) allow us to limit ourselves to check that
$v_P(\pi^{-1}(H))\ge q+1$ for every point $P$ chosen such that
$\pi(P)\in\P^N(\fq)$. As we have already noted, $v_P(f_i)\ge 1$ for $i\ge 1$.
Then, taking also into consideration Eq.~(\ref{eq2.3}), we only need to see
that $a^{(0)}_1=0$. As a matter of fact, this follows from Eq.~(\ref{eq2.2}),
and hence the proof of the lemma is complete.
        \end{proof}
As a corollary to the Lemmas \ref{lemma2.1} and \ref{lemma2.4} we obtain
the following result.

    \begin{theorem}\label{thm2.1} 
The morphism $\pi$ is a closed embedding, i.e., $\cX$ is $\fq$-isomorphic
to $\pi(\cX).$
    \end{theorem}

\begin{remark}\label{rem2.1} (1) As it was showed in \cite[Sect. 2]{FGT},
\cite[Sect. 2.3]{FT1}, a class of $\fq$-maximal curves can be characterized
by the type of the Weierstrass semigroup at some $\fq$-rational point of
the curve. The semigroups involved in such a characterization belong to a
special family of numerical semigroups $H$ defined by the following two
properties: (1) $q, q+1\in H$ and (2) there
exist $r,s\in H$ so that each $h\in H$ with $h\le q+1$ is generated by $r$
and $s$. Indeed, if a $\fq$-maximal curve has a $\fq$-rational point $P_0$
such that the Weierstrass semigroup $H(P_0)$ at $P_0$ satisfies each of the
above two conditions, then $H(P_0)=\langle r,s\rangle$. In particular, the
genus of such a curve is $(r-1)(s-1)/2$. Other interesting properties
of maximal curves depending on the behavior of their Weierstrass points
were pointed out in \cite[Sect. 2.4]{FT1}.

(2) Theorem \ref{thm2.1} implies that
    \begin{equation}\label{aut}
{\rm Aut}_{\fq}(\cX)\cong \{A\in PGL(N+1,q^2): A\pi(\cX)=\pi(\cX)\}\, .
    \end{equation}
For a stronger result on $\rm Aut(\cX)$, see Theorem \ref{thm3.3}.
     \end{remark}

\begin{remark}\label{rem2.2} For an application of Theorem \ref{thm2.1} in
Section \ref{s3} we stress that the condition of $\cD$ being a complete
linear series was not used. Hence Theorem \ref{thm2.1} holds true if $\cD$
is replaced by a (non-complete) linear subseries $\cR$ of $\cD$ as long as
$\cR$ contains all divisors $qP+\frx(P)$ with $P\in \cX$, and $\pi$ means
the morphism associated to $\cR$.
\end{remark}

\section{On the dual of $\pi(\cX)$}\label{s3}

The {\em dual} curve (also called {\em strict dual}) $\cZ^\star$ of a
non-degenerate projective geometrically irreducible algebraic curve $\cZ$
of a projective space $\P$ is the closure in the dual projective space
$\P^{\star}$ of the subset of points which represent the osculating
hyperplane $L_P^{N-1}$ to $\cZ$ at some general point $P\in\cZ$, see for
instance \cite{HK} and \cite{GV3}.

In this section, we assume that $\cX$ is a maximal curve over $\fq$, and
identify $\cX$ with $\pi(\cX)$ according to Theorem \ref{thm2.1}.  Let
$\pi^\star: \cX \rightarrow \P^N(\bar\fq)$ be the morphism with coordinate
functions $z_0,z_1\ldots,z_N$ introduced in the previous section. By
Lemma \ref{lemma2.3}(1), $\pi^\star\circ\frx$ is the Gauss map $P \mapsto
L_P^{(N-1)}$. This leads us to consider the curve $\pi^\star(\cX)$ in
$\P^N(\bar\fq)$. Note that, $\pi^\star(\cX)$ might be a degenerate curve
in the sense that it might happen that $\pi^\star(\cX)$ is contained and
non-degenerated in a subspace $\P^M$ of $ \P^N(\bar\fq)$.  By a result due
to Kaji \cite[Prop. 1]{K}, if this is the case then there is a 
$(N-M)$-dimensional subspace $\P^{N-M}$ of $\P^N(\bar\fq)$ which is the
intersection of the osculating hyperplane to $\cX$ at general points $P\in
\cX$, that is apart from a finite number of points $P\in \cX$.

In our situation, no point of $\cX$ lies on $\P^{N-M}$.  In fact, let
$R\in \P^{N-M}$, and assume on the contrary that $R\in \cX$.  Choose a
point $Q\in \cX$ such that $Q\neq R$ but the osculating hyperplane $L_Q$
to $\cX$ at $Q$ contains $P^{M-N}$. Since $L_Q$ meets $\cX$ in $\{Q,
\frx(Q)\}$ we have that $\frx(Q)=R$, and hence $Q$ is uniquely determined
by $R$. But this is a contradiction, as we can choose $Q$ in infinite
different ways.

Furthermore, $\P^M$ is invariant under the Frobenius collineation
$(X_0:\ldots :X_N)\mapsto (X^{q^2}_0:\ldots :X^{q^2}_N)$. This yields that
$\P^M$ is defined over $\fq$. Take a new $\fq$-invariant frame in
$\P^N(\bar\fq)$ in such a way that $\P^M$ has equation
$X_{M+1}=0,\ldots\,X_N=0$. Then $z_{M+1}=0,\ldots ,z_N=0$ and
$\pi^\star:\cX\to \P^M$ is given by $(z_0:\ldots:z_M)$. Hence, according
to Lemma \ref{lemma2.3}(1), the equation of the osculating hyperplane to
$\cX$ at $Q$ is $\gamma_0^q X_0+\ldots +\gamma_M^q X_M=0$, where
$\pi^\star(Q)=(\gamma_0:\ldots:\gamma_M)$.

        \begin{lemma}\label{lemma3.1} We have
$\deg(\pi^\star(\cX))=q+1$, and the linear series cut out on
$\pi^\star(\cX)$ by hyperplanes of $\P^M$ contains all divisors
$qP+\frx(P)$ with $P\in\cX$.
        \end{lemma}

\begin{proof} Choose a point $P_0=(\alpha_0:\ldots:\alpha_N)\in\cX(\fq)$.
Here, $\alpha_i\neq 0$ for some $i$ with $0\le i\le M$. In fact, if
$\alpha_i=0$ for $i=0,\ldots,M$, then $P_0$ would belong to the hyperplane
osculating at general points of $\cX$ and so $P_0$ would be in
the above $\P^{N-M}$ which is impossible as we have shown before. Now
consider the hyperplane $H$ of equation
$\alpha_0^q X_0+\ldots+\alpha_M^q X_M=0$ which can be regarded as a
hyperplane of
$\P^M$. Let $P\in\cX$ such that
$\pi^\star(P)=(\gamma_0:\ldots:\gamma_M)\in H\cap\pi^\star(\cX)$.
$\gamma_i\in\bar\fq$. We have that
$\alpha_0^q\gamma_0+\ldots+\alpha_M^q\gamma_M=0$ so that
$\gamma^q_0\alpha_0+\ldots+\gamma^q_M\alpha_M=0$.
This shows that the osculating hyperplane to $\cX$ at $P$ passes through
$P_0$ (Lemma \ref{lemma2.3}(1)). Since $P_0\in\cX(\fq)$, this is only
possible when $P=P_0$.
Thus we have proved that $H\cap\pi^\star(\cX)$ contains no point
different from $\pi^\star(P_0)$. We want to show next that the divisor
$(\pi^\star)^{-1}(H)$ of $\cX$ is $(q+1)P_0$. To do this we have
to show that
$$
v_{P_0}((\pi^\star)^{-1}(H))=v_{P_0}(\alpha^q_0 w_0+\ldots
+\alpha^q_Nw_N)=q+1\, ,
$$
where $v_{P_0}$ denotes the valuation at $P_0$, $w_i:=t^{e_{P_0}}z_i$, $t$
is a local parameter at $P_0$ and $e_{P_0}:=-{\rm
min}\{v_{P_0}(z_0),\ldots, v_{P_0}(z_N)\}$. (Recall that
$z_{M+1}=\ldots=z_N=0$.)

After a $\fq$-linear transformation of $\P^N(\bar\fq)$ we may assume that
$P_0=(1:0:\ldots:0)$ and that 
  $$
f_0=1,\quad f_1=t^{j_1}+\ldots,\quad f_N=t^{j_N}+\ldots\, ,
  $$
where $f_i=f_i/f_0$ and $(0,j_1,\ldots,j_N)$ is the $(\cD,P_0)$-order
sequence of $\cX$. Then we have to show that $v_{P_0}(w_0)=q+1$.

From Eq. (\ref{eq2.2}) we deduce that
  \begin{equation*}
w_0(t)+w_1(t)(t^{j_1}+\ldots)^q+\ldots+w_N(t)(t^{j_N}+
\ldots)^q=0\, .\tag{$*$}
  \end{equation*}

On the other hand, we claim that $v_{P_0}(w_1(t))=1$. By (\ref{eq2.1})
  $$
w_0(t)^q+w_1(t)^q(t^{j_1}+\ldots )+\ldots+w_N(t)^q(t^{j_N}+\ldots)=0\, .
   $$
From the definition of $w_i$ it follows that $v_{P_0}(w_i(t))=0$ for
almost one index $i$. Since $1=j_1<j_2<\ldots<j_N=q+1$ and $j_{N-1}\le q$
the only possibility is $i=N$, and $w_1(t)=ut+\ldots$ with $u\neq 0$.
The latter relation proves the claim. Now, this together with Eq. $(*)$
yield that $v_{P_0}((\pi^\star)^{-1}(H))=q+1$. Hence,
$(\pi^\star)^{-1}(H)$ of $\cX$ is $(q+1)P_0$ from which the first part of
the Lemma \ref{lemma3.1} follows. The second part follows from the
Fundamental Equivalence (\ref{eq1.2}).
\end{proof}

The above Lemma together with Remark \ref{rem2.2} have the following
corollary.
  \begin{lemma}\label{lemma3.2} The curves $\cX$ and $\pi^\star(\cX)$ are
$\fq$-isomorphic.
  \end{lemma}

Also, since $\cD$ is a complete linear series, Lemma \ref{lemma3.1}
gives the following result:
   \begin{lemma}\label{lemma3.3}
Every $z_i, 0\leq i\leq N$, is an $\fq$-linear combination of
$f_0,\ldots,f_N$.
   \end{lemma}
Now, we are in a position to prove the following theorem.

   \begin{theorem}\label{thm3.1} The curve
$\cX$ lies on a Hermitian variety defined over $\fq$ of $\P^N(\bar\fq)$.
   \end{theorem}

\begin{proof} Without loss of generality we may suppose that $f_0=z_0=1$.
For $i=0,\ldots, N$, let $z_i=\sum_{j=0}^{N}c_{ij}f_j$ with $c_{ij}\in
\fq$. Note that $c_{ij}=0$ for $M+1\leq i\leq N$ and that
the matrix $C=(c_{ij})$ has rank $M+1$. We prove that $C$ is
actually a Hermitian matrix over $\fq$. To do this, we re-write
Eq. (\ref{eq2.2}) in the following manner:
 $$
1+(\sum_{i=0}^{N}
(c_{i1}^qf_i)^qf_1+\ldots+
(\sum_{i=0}^{N}(c_{iN}^qf_i)^qf_N=0\, .
$$

Taking into account the uniqueness of the $N$-tuple
$(z_0=1,z_1,\ldots,z_N)$ proved in Lemma \ref{lemma2.3}(3), comparison
with Eq. (\ref{eq2.1}) gives

$$
\sum_{i=0}^{N} c_{i1}^qf_i=\sum_{i=0}^{N} c_{1i}f_i,\ldots,
\sum_{i=0}^{N} c_{iN}^qf_i=\sum_{i=0}^{N} c_{Ni}f_i\, .
$$

Since $f_0=1,f_1,\ldots,f_N$ are linearly independent over $\fq$, this
yields $c_{ij}=c_{ji}^q$ for every $0\leq i,k\leq N$. This proves that $C$
is Hermitian. After a $\fq$-linear transformation of $\P^N(\bar\fq)$ we
may assume that the matrix $C$ is the diagonal matrix with $M$ units. Then
(\ref{eq2.1}) becomes
        $$
f^{q+1}_0+\ldots+f^{q+1}_M=0\, ,
        $$
and hence $\cX$ lies on the Hermitian variety of equation
$X^{q+1}_0+\ldots+X^{q+1}_M=0.$
   \end{proof}
       \begin{remark}\label{rem3.1}
From the proof above, $z_i=f_i$ for $0\leq i\leq M$. Hence
$\pi^\star(\cX)$ is the projection $(f_0:\ldots:f_N) \rightarrow
(f_0:\ldots:f_M)$, and $\pi^\star(\cX)$ lies on a
non-degenerate Hermitian variety defined over $\fq$ of $\P^M(\bar\fq)$.
        \end{remark}
Taking into account Lemma \ref{lemma3.2} we obtain the following result.
        \begin{theorem}\label{thm3.2}
$\cX$ admits a non-singular model given by a curve defined over $\fq$
which has degree $q+1$ and lies on a non-degenerate Hermitian variety
defined over $\fq$ of $\P^M(\bar\fq)$ of dimension $M\leq N$.
        \end{theorem}

From the above arguments, it also turns out that the osculating hyperplane
to $\cX$ at any point $P\in\cX$ coincides with the tangent hyperplane to
the non-degenerate Hermitian variety at the same point $P$. This allows us
to improve the previous result (\ref{aut}) on ${\rm Aut}_{\fq}(\cX)$:

        \begin{theorem}\label{thm3.3}
${\rm Aut}_{\fq}(\cX)$ is isomorphic to a subgroup of the projective
unitary group $PGU(M+1,\fq)$.
        \end{theorem}
        \begin{proof}
By a way of contradiction, assume that $\cX$ lies not only on $\cH$ but also
on the non-degenerate Hermitian variety $\cH^{\prime}$ which is assumed to be
the image of $\cH$ by a non-trivial $\fq$-linear collineation fixing $\cX$.
Choose any point $P\in\cX$. Then the $\cH$ and $\cH^{\prime}$ have the same
tangent hyperplane at $P$, as each of these tangent hyperplanes coincides with the osculating
hyperplane to $\cX$ at $P$. To express this geometric condition
in algebraic terms, set $P:=(x_0:\ldots:x_M)$, and write
the equations of $\cH$ and $\cH^{\prime}$ explicitly:
$\cH:=X^{q+1}_0+\ldots+X^{q+1}_N=0$; $\cH^{\prime}:={\bf X}^tC({\bf X})^q=0$
where ${\bf X}:=(X_0,\ldots,X_M)$, and $C$ is a non-singular non-identity
unitary matrix of rank $M+1$. Then the above geometric condition in
algebraic terms is that the homogeneous $(M+1)$-tuples
$(x_0^q,\ldots,x_M^q)$ and
$(c_{0,0}x^q_0+\ldots+c_{M,0}x^q_M,\ldots,c_{0,M}x^q_0+\ldots+c_{M,M}x^q_M)$
are equal up to a non-zero factor. Another meaning of the latter relation is
that the non-trivial $\fq$-linear collineation associated to the matrix $C$
fixes $\cX$ pointwise. But this is impossible as $\cX$ is not contained in a
hyperplane of $\bf P^M$; a contradiction which proves the theorem.
        \end{proof}    
  
        \section{Curves lying on a Hermitian variety}\label{s4} The aim of
this section is to show that the property given in Theorem \ref{thm3.2}
characterizes $\fq$-maximal curves. For this purpose,
we assume from now on that $\cX$ is a projective geometrically irreducible
non-singular algebraic curve defined over a finite field $\fq$ which is
equipped with a non-degenerated $\fq$-birational morphism
$\pi=(f_0:\ldots:f_M): \cX\to \P^M(\bar\fq)$ such that the curve
$\cY:=\pi(\cX)$ has the following properties:
    \begin{itemize}
\item It has degree $q+1$, and it lies on a non degenerate Hermitian
variety $\cH\subseteq \P^M(\bar\fq)$ defined over $\fq$.
    \end{itemize}
The main result in this section is the following theorem.

        \begin{theorem}\label{thm4.1}
The curve $\cX$ is $\fq$-maximal.
        \end{theorem}

The Hermitian variety $\cH$ of $\P^M(\bar\fq)$ is assumed to
be in its canonical form $X^{q+1}_0+\ldots+X^{q+1}_M=0$. By our hypothesis,

        \begin{equation}\label{eq4.1}
f^{q+1}_0+\ldots+f^{q+1}_M=0\, .
        \end{equation}

For any point $P\in\cX$, let $\pi(P)=(\alpha_0:\ldots:\alpha_M)$.
Choose a local parameter $t$ at $P$, and arrange the coordinate functions to
have $v_P(f_i)\geq 0$ for $i=0,\ldots,M$ and
$v_P(f_k)=0$ for at least one index $k\in\{0,\ldots,M\}$. Then
   $$
f_i(t)= \sum_{j=0}^{\infty}a_{i,j}t^j\in \bar\fq[[t]]
   $$
is the local expansion of $f_i$ at $P$. Here,
$\alpha_i=a_{i,0}$ and $a_{k,0}\neq 0$. The tangent hyperplane $H_P$ to the
Hermitian variety at $\pi(P)$ has equation
$\alpha^q_0X_0+\ldots+\alpha^q_MX_M=0$.

The first step toward Theorem \ref{thm4.1} is the following lemma.

        \begin{lemma}\label{lemma4.1}
The linear series $\cR$ cut out on $\cY$ by hyperplanes contains the
divisor $qP+\frx(P)$ for every $P\in\cX$.
        \end{lemma}
        \begin{proof} We show that $H_P$ cuts out on $\cY$ the divisor
$qP+\frx(P)$. From Eq. (\ref{eq4.1}),

        \begin{equation}\label{eq4.2}
(\sum_{j=0}^{\infty}a_{0,j}t^j)^qf_0+
\ldots+(\sum_{j=0}^{\infty}a_{M,j}t^j)^qf_M=0\, .
        \end{equation}

Writing the lower order terms in $t$, we have
    $$
\sum_{i=0}^{M}a_{i,0}^qf_i+t^q\sum_{i=0}^{M}a_{i,0}a_{i,1}^q+
t^{q+1}\sum_{i=0}^{R}a_{i,1}^{q+1}+t^{q+2}[\ldots]=0\, .
    $$
Hence $v_P(\pi^{-1}(H_P)\ge q$
and equality holds if and only if $\sum_{i=0}^{M}a_{i,1}^qa_{i,0}\neq0$.
We show that if $P\in\cX(\fq)$, then $\sum_{i=0}^{M}a_{i,1}^qa_{i,0}=0$. 
From (\ref{eq4.2}),
   $$
\sum_{i=0}^{M}a_{i,0}^{q+1}+t\sum_{j=0}^{M}a_{i,0}^qa_{i,1}+
t^q[\ldots]=0\, .
   $$
Thus, $\sum_{i=0}^{M}a_{i,0}^qa_{i,1}=0$. Since
$(\sum_{i=0}^{M}a_{i,0}^qa_{i,1})^q=\sum_{i=0}^{M}a_{i,0}a_{i,1}^q$
for $P\in\cX(\fq)$, the claim follows. Since $\pi$ is birational and
$\deg(\cY)=q+1$, we obtain $\pi^{-1}(H_P)=(q+1)P$ for every
$P\in\cX(\fq)$, which shows the lemma for every $P\in\cX(\fq)$.
For the case $P\not\in\cX(\fq)$, we also need to check that $\frx(P)\in H_P$.
This inclusion occurs when $\sum_{i=0}^{M}\alpha^{q^2+q}_i=0$.
Since the latter relation is a consequence of (\ref{eq4.2}), the claim
follows. Hence, $\pi^{-1}(H_P)=qP+\frx(P)$ because $\pi$ is birational and
$\deg(\cY)=q+1$,
        \end{proof}
Then from Remark \ref{rem2.2} and Lemma \ref{lemma4.1} follow that $\cX$
and $\cY=\pi(\cX)$ are $\fq$-isomorphic. Hence if $M=2$, $\cY$ is the
Hermitian curve and so $\cX$ is $\fq$-maximal. From now on we assume $M\ge
3$.

Our approach is based on a certain relationship between the Wronskians
determinants of $\cY$ and of its projection on a $(M-1)$-dimensional
subspace of $\P^M(\bar\fq)$. More precisely, let
$\bar\pi:\cX\to\P^{M-1}(\bar\fq)$ defined by $\cX\to
(f_0:\ldots :f_{M-1})$,
that is $\bar\cY$ is the projection of $\cY$ from the point
$(0:\ldots0:1)$ on the hyperplane $X_M=0$. It might happen that $\cY$ and
$\bar\cY$ are not $\fq$-birationally equivalent. However, it is always
possible avoid this situation by changing the coordinate system in
$\P^M(\bar\fq)$, see Appendix. So we assume that $\cY$ is
$\fq$-birationally equivalent to $\bar\cY$.

Choose a separating variable $t$ of $\cX$, and define $D_t$ as the Hasse
derivative with respect to $t$, see \cite{He}. Then \cite[Section 1]{SV}:

   $$
Wr(f_0,\ldots, f_{M-1}):=\det \left[\begin{array}{cccc}
D_t ^{\gep_0}f_0 & D_t ^{\gep_0}f_1 &\ldots & D_t ^{\gep_0}f_{M-1}
\\
\vdots & \vdots & \ldots & \vdots\\
D_t ^{\gep_{M-1}}f_0 & D_t ^{\gep_{M-1}}f_1 &\ldots & D_t
^{\gep_{M-1}}f_{M-1}
\end{array}\right]\, ,
   $$
and
   $$
Wr(f_0,\ldots, f_M):=det \left[\begin{array}{cccc}
D_t ^{\gep_0}f_0 & D_t ^{\gep_0}f_1 &\ldots & D_t
^{\gep_0}f_M
\\
\vdots & \vdots & \ldots & \vdots\\
D_t ^{\gep_{M}}f_0 & D_t ^{\gep_{M}}f_1 &\ldots &
D_t ^{\gep_{M}}f_M
\end{array}\right]\, .
    $$
Note that in our case $\gep_0=0,\gep_1=1,\ldots,\gep_M=q$.
        \begin{lemma}\label{lemma4.2} We have that
   $$
\div(Wr(f_0,\ldots,f_M))=\div(Wr(f_0,\ldots,f_{M-1}))-q\div(f_M)+
\div(f_0D_t^qf_0^q+\ldots + f_MD_t^qf_M^q)\, .
   $$
        \end{lemma}
    \begin{proof} Multiplying the last column by $f_M^q$ and adding to it
$f_0^q$ times the first column plus $f_1^q$ times the second column
etc. plus $f_{M-1}^q$ times the penultimate column gives
   $$
f_M^qWr(f_0,\ldots\,f_M)=
\left[\begin{array}{cccc}
f_0 & f_1 &\ldots &f_0^{q+1}+\ldots+f_M^{q+1} \\
D_t f_0 & D_t f_1 &\ldots &f_0^qD_tf_0+\ldots+f_M^qD_tf_M
\\
\vdots & \vdots & \ldots & \vdots\\
D_t ^{q}f_0 & D_t ^{q}f_1 &\ldots
&f_0^qD_t^qf_0+\ldots+f_M^qD_t^qf_M
\end{array}\right]\, .
    $$
Each element but the last one in the last column is actually $0$. In fact,
this follows from the relation (\ref{eq4.1}) by derivation. Furthermore,
the $q$-th Hasse derivative of the same relation gives
     $$
f_0^qD_tf_0+\ldots+f_M^qD_tf_M+f_0D_t^qf_0^q+\ldots+f_MD_t^qf_M^q=0\, ,
     $$
and this completes the proof.
     \end{proof}
Let $R_M$ be the ramification divisor of the linear series cut out
on $\cY$ by hyperplanes of $\P^M(\bar\fq)$. The following result
comes from \cite[p.6]{SV}:
      \begin{lemma}\label{lemma4.3}
Let $P\in\cX$. If $t$ is a local parameter of $\cX$ at $P$, then
    $$
v_P(R_M)=v_P(Wr(f_0,\ldots,f_M))\, .
    $$
       \end{lemma}
Similarly, let $R_{M-1}$ be the ramification divisor of the linear
series cut out on $\bar\cY$ by hyperplanes of $\P^{M-1}(\bar\fq)$.
      \begin{lemma}\label{lemma4.4}
Let $P\in\cX$. If $t$ is a local parameter of $\cX$ at $P$, then
   $$
v_P(R_{M-1})=v_P(Wr(f_0,\ldots,f_{M-1}))\, .
    $$
     \end{lemma}
     \begin{proof} By \cite[p.6]{SV},
  $$
v_P(R_{M-1})=v_P(Wr(f_0,\ldots,f_{M-1}))+
(\gep_0+\gep_1+\ldots\gep_{M-1})v_P(dt)+M\bar e_P\, ,
  $$
where $\bar e_P:=-\min\{v_P(f_0),\ldots ,v_P(f_{M-1})\}$. Actually,
$\bar e_P=0$. In fact, $\bar e_P>0$ together with $e_P=0$ would imply
that the point $U_M:=(0:\ldots :0:1)$ lies on $\cY$ but this
contradicts (\ref{eq4.1}). Since $t$ is a local parameter at $P$, we also
have $v_P(dt)=0$, and the claim follows.
      \end{proof}
The following result will play a crucial role in the sequel.
     \begin{lemma}\label{lemma4.5}
$$
v_P(f_0D_t^qf_0^q+\ldots+f_MD_t^qf_M^q)
=\left\{ \begin{array}{ll}
  1 & \mbox{when $P\in\cX(\fq)$\, ,}\\
  0 & \mbox{when $P\not\in\cX(\fq)$\, .}
 \end{array}\right.
$$
       \end{lemma}
       \begin{proof} From the proof of Lemma \ref{lemma4.1} we obtain the
following result. For any point $P\in\cX$,
        \begin{itemize}

\item $P\not\in\cX(\fq)$ if and only if $\sum_{i=0}^{M}a_{i,1}^qa_{M,1}\neq0,$
\item $P\in\cX(\fq)$ if and only if $\sum_{i=0}^{M}a_{i,1}^qa_{M,1}=0$ but
$\sum_{i=0}^{M}a_{i,1}^{q+1}\neq 0.$
        \end{itemize}
On the other hand,
   $$
f_0D_t^qf_0^q+\ldots+f_MD_t^qf_M^q=
(\sum_{j=0}^{\infty}a_{0,j}t^j)(a_{0,1}^q+t^q[\ldots])+\ldots+
(\sum_{j=0}^{\infty}a_{M,j}t^j)(a_{M,1}^q+t^q[\ldots])\, .
   $$
Hence
        \begin{itemize}
\item $v_P(f_0D_t^qf_0^q+\ldots+f_MD_t^qf_M^q)=0$ if and
only if $\sum a_{i,0}a_{i,1}^q\neq0.$
\item $v_P(f_0D_t^qf_0^q+\ldots+f_MD_t^qf_M^q)=1$ if and only if
$\sum a_{i,0}a_{i,1}^q=0$ but $\sum a_{i,1}^{q+1}\neq0.$
        \end{itemize}
Now, comparison with the previous result proves Lemma \ref{lemma4.5}.
     \end{proof}
%
Now we are in a position to finish the proof of Theorem \ref{thm4.1}. By
\cite [p.6]{SV},
   $$
\sum v_P(R_M)=(\gep_0+\gep_1+\ldots\gep_{M})(2g-2)+(M+1)(q+1)\, ,
   $$
and
   $$
\sum v_P(R_{M-1})=(\gep_0+\gep_1+\ldots\gep_{M-1})(2g-2)+M(q+1)\, .
   $$
Hence $\sum (v_P(R_M)-v_P(R_{M-1}))=q(2g-2)+q+1$. Lemmas
\ref{lemma4.2}, \ref{lemma4.3}, \ref{lemma4.4}, and
\ref{lemma4.5} together with $\sum v_P(f_M)=q+1$ give Theorem
\ref{thm4.1}.
        \section{Examples}\label{s5}
We will show how each of the known examples of maximal curves with $\deg(
\cD)=3$ can be embedded in a non degenerate Hermitian variety defined
over $\fq$ of
${\P^3(\fq)}$. In this way we obtain an independent proof of the
maximality of these curves.
   \begin{example}\label{ex5.1} (\cite[Thm. 2.1.(IV)(2)]{CKT1}) Let
$q\equiv 2\pmod{3}$, and fix a primitive third root of unity
$\epsilon\in\fq$. For $i=0,1,2$, let $\cC_i$ be a projective geometrically
irreducible non-singular curve defined over $\fq$ whose
function field over $\fq$ is generated by $x$ and $y$ satisfying the
irreducible polynomial relation
    $$
        \label{eqes1}
\epsilon^i x^{(q+1)/3}+\epsilon^{2i} x^{2(q+1)/3}+y^{q+1}=0\, .
    $$
Let
     $$
f_0:=x;f_1:=x^2;f_2:=y^3,f_3:=xy
     $$
be the coordinate functions of a morphism $\pi=\cC_i \rightarrow
\P^3(\fq)$. Note that the these three curves $\pi(\cC_i)$
are pairwise projectively equivalent in $\P^3(\fq)$. In fact, the
linear transformation induced by the matrix
    $$
T^{(i)}_4= \begin{bmatrix}
\epsilon^i & 0 & 0 & 0 \\
0 & \epsilon^{2i} & 0 & 0 \\
0 & 0 & 1 & 0 \\
0 & 0 & 0 & \epsilon^i
\end{bmatrix}\,
    $$
maps $\pi(\cC_0)$ to $\pi(\cC_i)$. We show that $\pi(\cC_i)$
is projectively equivalent to a projectively irreducible non-singular curve defined over $\fq$
and contained in the non-degenerate Hermitian surface $\cH_3$ of equation
$X^{q+1}_0+X^{q+1}_1+X^{q+1}_2+X^{q+1}_3=0$.
To do this we start with the relation in $\fq[X,Y]$
        \begin{align*}
        \begin{split}
(X^{(q+1)/3}+X^{2(q+1)/3}+Y^{q+1})(\epsilon X^{(q+1)/3}+
\epsilon^2 X^{2(q+1)/3}+Y^{q+1})\\
(\epsilon^2 X^{(q+1)/3}+\epsilon X^{2(q+1)/3}+Y^{q+1})
=X^{q+1}+X^{2(q+1)}+Y^{3(q+1)}-3X^{q+1}Y^{q+1}
       \end{split}
        \end{align*}
which yields $x^{q+1}+x^{2(q+1)}+y^{3(q+1)}-3x^{q+1}y^{q+1}=0$. Thus,
$f_0^{q+1}+f_1^{q+1}+f_2^{q+1}-3f_3^{q+1}=0$. This shows that
$\pi(\cC_i)$ lies on the non degenerate Hermitian variety $\cH$ of equation
$X_0^{q+1}+X_1^{q+1}+X_2^{q+1}+X_3^{q+1}=0$, up to the linear transformation
$(X_0,X_1,X_2,X_3)\rightarrow (X_0,X_1,X_2,wX_3)$ with $w^{q+1}=-3$.
Furthermore, $\pi(\cC_i)$ is contained in the cubic surface $\Sigma_3$
of $P^3(\fq)$ of equation $X_3^3+w^3X_0X_1X_2=0$. More precisely,
the intersection curve of $\cH$ and $\Sigma_3$ splits
into the above three pairwise projectively equivalent curves,
namely $\pi(\cC_0)$, $\pi(\cC_1)$, and $\pi(\cC_2)$, each of degree $q+1$.
By Theorem \ref{thmC}, $\pi(\cC_i)$ is a non-singular maximal curve
defined over $\fq$. According to \cite[Thm. 2.1.(IV)(2)]{CKT1},
its genus is equal to $(q^2-q+4)/6$.
    \end{example}
    \begin{example}\label{ex5.2} (\cite[sect. 6]{CKT}) A similar but non
isomorphic example is given in \cite{CKT}.
Again, assume that $3$ divides $q+1$, and fix a primitive third
root of unity $\epsilon\in\fq$. For $i=0,1,2$, let $\cC_i$ be curves as in
Example \ref{ex5.1} whose function field over $\fq$ is generated by $x$
and $y$ satisfying the irreducible polynomial relation
     $$
        \label{eqes2}
\epsilon^i yx^{(q-2)/3}+y^q+\epsilon^{2i} x^{(2q-1)/3}=0\, .
     $$
Let
    $$
f_0:=x;f_1:=x^2;f_2:=y^3,f_3:=-3xy
    $$
be the coordinate functions of a morphism $\pi=\cC_i \rightarrow
\P^3(\fq)$. Note that the these three curves $\pi(\cC_i)$
are pairwise projectively equivalent in $\P^3(\fq)$. In fact, the
linear transformation induced by the matrix
    $$
T^{(i)}_4= \begin{bmatrix}
\epsilon^i & 0 & 0 & 0 \\
0 & \epsilon^{2i} & 0 & 0 \\
0 & 0 & 1 & 0 \\
0 & 0 & 0 & \epsilon^i
\end{bmatrix}\,
    $$
maps $\pi(\cC_0)$ to $\pi(\cC_i)$.

We show that $\pi(\cC_i)$
is projectively equivalent to a projectively irreducible non-singular
curve defined over $\fq$
and contained in the non-degenerate Hermitian surface $\cH_3$ of equation
$X^{q+1}_0+X^{q+1}_1+X^{q+1}_2+X^{q+1}_3=0$.
To do this we start with the relation in $\fq[X,Y]$
        \begin{align*}
        \begin{split}
(YX^{q-3}+Y^q+X^{(2q-1)/3})(\epsilon YX^{q-3}+Y^q+\epsilon^2 X^{(2q-1)/3})
(\epsilon^2 YX^{q-3}+Y^q+\epsilon X^{(2q-1)/3})\\
=Y^3X^{q-2}+Y^{3q}+X^{2q-1}-3X^{q-1}Y^{q+1}
       \end{split}
        \end{align*}
which implies $y^3x^{q-2}+y^{3q}+x^{2q-1}-3x^{q-1}y^{q+1}=0$. Thus,
$y^3x^q+y^{3q+2}+x^{2q+1}-3x^{q+1}y^{q+1}=0$, and hence
$f_2f^q_0+f^q_2f_1+f^q_1f_0-3f^{q+1}_3=0$. This shows that
$\pi(\cC_i)$ lies on the surface $\Sigma_{q+1}$ of equation
$X^q_0X_1+X^q_1X_2+X^q_2X_0-3X^{q+1}_3=0$. Furthermore,
$\pi(\cC_i)$ is contained in the cubic surface $\Sigma_3$
of $P^3(\fq)$ of equation $X_3^3+27X_0X_1X_2=0$. More precisely,
the intersection curve of $\Sigma_{q+1}$ and $\Sigma_3$ splits
into the above three pairwise projectively equivalent curves,
namely $\pi(\cC_0)$, $\pi(\cC_1)$, and $\pi(\cC_2)$, each of degree $q+1$.

To prove that $\Sigma_{q+1}$ is projectively equivalent to $\cH_3$,
we use the same argument employed in \cite {CKT}.
Choose a root $a$ of the polynomial $p(X):=X^{q+1}+X+1$. Then
$a^{q^2+q+1}=1$, and hence $a\in\mathbf F_{q^3}$. By \cite[Lemma 4]{CKT},
$a^{q+1}+a^{q^2+q+1}+a=0$ and $a^{q^2+q+2}+a^{q+1}+1=0$,
but $a^{q+2}+a^{q^2+1}+a^q\neq 0$ as
$(a^{q+2}+a^{q^2+1}+a^q)^{q-1}=a^{-1}$. Furthermore, the matrix
    $$
M_3= \begin{bmatrix}
a & 1 & a^{q^2+1}  \\
a^{q^2+1} & a & 1  \\
1 & a^{q^2+1} & a  \\
\end{bmatrix}\,
    $$
is non-singular. Also, choose an element $\mu\in\mathbf F_{q}$ satisfying
$-3\mu^{q+1}=a^{q^3+q+1}+a^{q^2+1}+a^q$, and define $\kappa$ as the projective linear
transformation $\kappa: P^3(\bar\mathbf F_q) \rightarrow
P^3(\bar\mathbf F_q)$ induced by the non--singular matrix
    $$
M_4= \begin{bmatrix}
a & 1 & a^{q^2+1} & 0 \\
a^{q^2+1} & a & 1 & 0 \\
1 & a^{q^2+1} & a & 0 \\
0 & 0 & 0 & -\mu \\
\end{bmatrix}\, .
    $$
A straightforward computation shows that $\kappa^{-1}$ maps $\Sigma_{q+1}$
to $\cH_3$, and $\Sigma_3$ to the cubic surface $\bar\Sigma_3$ of equation
    $$
(X^3_0+X^3_1+X^3_2)+Tr[a^{q+1}](X^2_0X_1+X^2_1X_2+X^2_2X_0)
+Tr[a](X^2_0X_2+X^2_1X_0+X^2_2X_1)$$
$$+(3+Tr[a^{q-1}])X_0X_1X_2-a^{q-1}\mu^3X^3_3=0
    $$
where $Tr[u]:=u+u^q+u^{q^2}$ is the trace of $u\in\mathbf
F_{q^3}$ Furthermore $a^{q-1}\mu^3\in\fq$, and this shows that $\bar\Sigma_3$
is actually defined   over $\fq$.
Now, $\pi(\cC_i)$ is mapped under $\kappa^{-1}$ to a projectively irreducible
algebraic curve of degree $q+1$ defined over $\fq$ and contained in $\cH_3$.
By Theorem \ref{thmC} $\kappa^{-1}(\cC_i)$ is a non-singular maximal curve
defined over $\fq$. By \cite[Lemma 6.1.(5)]{CKT}, its genus is
equal to $(q^2-q-2)/6$.
    \end{example}
    \begin{example}\label{ex5.3} (\cite{FGT}) Let $q$ be odd
and for $i=1,2$, let
$\cC_i(\fq)$ be curves as in Example \ref{ex5.1} whose function field over
$\fq$ is generated by $x$ and $y$ such that
    $$
y^q+y+(-1)^ix^{(q+1)/2}=0\, .
    $$
The functions
    $$
f_0:=1;f_1:=x;f_2:=y;f_3:=y^2
    $$
define a morphism
$\pi=\cC_i \rightarrow\P^3(\fq)$. The resulting curves
$\pi(\cC_i)$ are projectively equivalent, since the linear
transformation induced by the matrix
    $$
T^{(i)}_4= \begin{bmatrix}
1 & 0 & 0 & 0 \\
0 & \epsilon & 0 & 0 \\
0 & 0 & 1 & 0 \\
0 & 0 & 0 & 1
\end{bmatrix}\,
    $$
with $\epsilon^{(q+1)/2}=-1$, maps $\pi_0(\cC_0)$ to $\pi(\cC_1)$.
The polynomial relation
$$(Y^q+Y-X^{(q+1)/2})(Y^q+Y+X^{(q+1)/2})=Y^{2q}+2Y^{q+1}+Y^2-X^{q+1}
$$
implies that $y^{2q}+2y^{q+1}+y^2-x^{q+1}=0$ in
$\cC_i(\fq)$. Hence $f_3^q+f_3+2f_2^{q+1}-f_1^{q+1}=0$ in
$\cC_i(\fq)$. This proves that $\pi(\cC_i)$ lies on the surface $\Sigma$ of
equation $X_3^qX_0+X_3X_0^q+2X_2^{q+1}-X_1^{q+1}=0$ which is
actually a non degenerate Hermitian variety defined over $\fq$ of 
$\P^3(\bar\fq)$. 
Also, $\cC_i$ lies on the quadratic cone $K$ of equation
$X_2^2-X_0X_3=0$, and hence the intersection of $\Sigma$ and $K$
splits into the curves $\pi(\cC_0)$ and $\pi(\cC_1)$. By Theorem
\ref{thmC} $\pi(\cC_i)$ is a a non-singular maximal curve
defined over $\fq$. Its genus is equal to $(q-1)^2/4$, according to \cite
{FGT}.
     \end{example}
     \begin{example}\label{ex5.4} (\cite{AT}) Let $q=2^t$, and put
$Tr[Y]:=Y+Y^2+\ldots+Y^{q/2}$. For $i=0,1\in\mathbf F_2\subset\fq$,
let $\cC_i$ be curves as in Example \ref{ex5.1} whose function field
over $\fq$ is generated by $x$ and $y$ such that
   $$
Tr[y]+x^{q+1}+i=0\, .
   $$
Let $\pi=\cC \rightarrow \P^3(\fq)$ be given by the coordinate
functions
   $$
f_0:=1;f_1:=x;f_2=y;f_3=x^2\, .
   $$
Since
   $$
(Tr[Y]+X^{q+1})+(Tr[Y]+X^{q+1}+1)=Y^q+Y+X^{q+1}+X^{2q+2}
   $$
we have
$y^q+y+x^{q+1}+x^{2q+2}=0$ in $\cC(\fq)$.
This implies that $\pi(\cC_i)$ lies on the non degenerate Hermitian variety
$H$ of equation $X_2^qX_0+X_2X_0^q+X_1^{q+1}+X_3^{q+1}=0$. Furthermore, the
quadratic cone $K$ of equation $X_3X_0=X_1^2$ also contains
$\pi(\cC_i)$. Hence $H\cap K$ splits into $\pi(\cC_0)$ and $\pi(\cC_1)$.
Note that $\pi(\cC_0)$ and $\pi(\cC_1)$ are projectively equivalent
curves in $\P^3(\fq)$, and hence both have degree $q+1$. Again by Theorem
\ref{thmC} $\pi(\cC_i)$ is a
non-singular maximal curve defined over $\fq$. Its genus is equal to
$q(q-2)/4$.
      \end{example}
      \begin{example}\label{ex5.5} (\cite{GeV}) Let $q=3^t$, and put
$Tr[Y]:=Y+Y^3+\ldots+Y^{q/3}$.
For $i=0,1,2\in\mathbf F_3\subset\fq$, let $\cC_i$ be curves as in Example
\ref{ex5.1} whose function field over $\fq$ is generated by $x$ and $y$
such
that
   $$
Tr[y]^2-x^q-x+i(Tr[y]+i)=0\, .
   $$
Since $(Tr[Y]^2+X^q-X)(Tr[Y]^2-X^q-X+Tr[Y]+1)(Tr[Y]^2-X^q-X-Tr[Y]+1)=
(X^q+X)(X^q+X-1)^2-(Y^q-Y)^2$, we have
    \begin{equation*}
(x^3+x^2-y^2+x)^q+(x^3+x^2-y^2+x)-x^{q+1}-y^{q+1}=0\, .\tag{$*$}
    \end{equation*}
Let $\pi=\cC_i \rightarrow \P^3(\fq)$ be given by the coordinate
functions $f_0:=1;f_1:=x;f_2=y,f_3:=x^3+x^2-y^2+x.$ It can be
checked that these three curves are pairwise projectively equivalent in
$\P^3(\fq)$. From Eq. ($*$), $\pi(\cC)$ lies on the non degenerate
Hermitian variety of equation $X_0X_3^q+X_0^qX_3-X_1^{q+1}-X_2^{q+1}=0$.
Furthermore, the cubic surface of equation $X_3X_0^2-X_1^3+X_1^2X_0+X_2^2X_0-X_1X_0^2$ also
contains $\pi(\cC_i)$. It turns out that $\pi(\cC_i)$ has degree $q+1$, and
Theorem \ref{thmC} ensures that $\pi(\cC)$ is a
non-singular maximal curve defined over $\fq$. Its genus is equal
to $q(q-1)/6$.
     \end{example}
     \begin{remark}\label{rem5.1} 
In all the above examples $\cX$ lies not only on a non-degenerate Hermitian
surface but also on a cubic surface. This is related to a classical result of
Halphen on reduced and irreducible complex algebraic curves in $\bf P^3$
not lying on a quadratic surface which states that the degree $d$ and the
genus $g$ satisfy of such a curve satisfy the following
inequality:

\begin{equation*}
g\leq\pi_1(d,3)= \left \{ \begin{array}{lll}
            d^2/6-d/2+1     & \mbox{\em for $d\equiv0$ mod $3;$\/} \\
            d^2/6-d/2+1/3 & \mbox{\em for $d\not\equiv0$ mod $3.$ \/}
                         \end{array}
          \right. \
\end{equation*}

A rigorous proof of the Halphen theorem and its extension to higher
dimensional spaces is found in the book \cite{EH}. Rathmann \cite{RJ}
pointed out that the proof also works in positive characteristic apart from some
possible exceptional cases related to the monodromy group of the
curve.
\end{remark}

     \section{Appendix}\label{s6}
For $M\ge 3$, in Section \ref{s4} we have claimed that the
curves $\cY\subseteq \P^M(\bar\fq)$ and $\bar\cY\subseteq
\P^{M-1}(\bar\fq)$ are
$\fq$-birationally equivalent, up to a change of coordinates in
$\P^M(\bar\fq)$. Now we give the proof. Notation and terminology being as
in Section
\ref{s4}, two technical lemmas are needed.
        \begin{lemma}\label{lemma6.1} The space
$\P^M(\fq)$ contains a point $P$ satisfying each of the following
three conditions:
        \begin{itemize}
\item $P$ is not on $\cH$;
\item no tangent line to $\cY$ at a $\fq$-rational point
passes through $P$;
\item no chord through two $\fq$-rational points of $\cY$
passes through $P$.
        \end{itemize}
        \end{lemma}
        \begin{proof} Take a $\fq$-rational point $Q\in\cY$. Since the
number of $\fq$-rational points of $\cX$ is $q^2+1+2gq\leq q^3+1$, there are
at most $q^3$ chords through $Q$ and another $\fq$-rational point
of $\cY$. But, since $M\ge 3$, the number of $\fq$-rational lines through
$Q$ is at least $q^4+q^2+1$ and hence one of these lines is neither a line
contained in $\cH$, nor a tangent line to $\cY$ at $Q$, nor a chord
through $Q$ and another $\fq$-rational point of $\cY$.
Now, any $\fq$-rational point $P$ outside $\cH$ is a good choice
for $P$.
        \end{proof}
        \begin{lemma}\label{lemma6.2} Let $r$ be a $\fq$-rational line
through a $\fq$-rational point $R$ of $\cY$. Then $r\cap\cY$ only
contains $\fq$-rational points from $\cY$.
        \end{lemma}
        \begin{proof}
Assume on the contrary that $r$ meets $\cY$ in a non $\fq$-rational
point $S$. Then $r$ is the line joining $S$ and ${\mathbf Fr}(S)$.
This implies that $r$ is contained in the osculating hyperplane of
$\cY$ at $S$. Hence the common points of $r$ with $\cY$
are only two, namely $S$ and ${\mathbf Fr}(S)$. But this
contradicts the hypothesis that $R\in r\cap\cY$.
        \end{proof}
Take a point $P$ as in Lemma \ref{lemma6.1}. By a classical result
(see \cite{SB}, and also \cite [23.4]{HT}), the linear collineation group
$PGU(M+1,q^2)$ preserving $\cH$ acts transitively on the set of
all points of $\P^M(\bar\fq)$ not on $\cH$. Hence a linear
collineation of $\P^M(\bar\fq)$ can be applied which preserves $\cH$ and
maps $P$ to $(0:\ldots:0:1)$. Lemmas \ref{lemma6.1} and \ref{lemma6.2}
ensure now that $\cY$ and $\bar\cY$ are $\fq$-birationally equivalent.


\begin{thebibliography}{99}

\bibitem{AT} M. Abd\'on and F. Torres, {\em On maximal curves in
characteristic two,} Manuscripta Math. {\bf 99} (1999), 39-53.


\bibitem{CKT} A. Cossidente, G. Korchm\'{a}ros and F. Torres, {\em
On curves covered by the Hermitian curve}, J. Algebra {\bf 216}
(1999), 56-76.

\bibitem{CKT1} A. Cossidente, G. Korchm\'{a}ros and F. Torres, {\em
Curves of large genus covered by the Hermitian curve}, to appear in Comm.
Algebra.

\bibitem{EH} D. Eisenbud and J. Harris, ``Curves in projective
space'', Les Presses de l'Universit\'e de Montr\'eal, Montr\'eal, 1982.

\bibitem{FGT} R. Fuhrmann, A. Garcia and F. Torres, {\em On maximal
curves}, J. Number Theory {\bf 67} (1997), 29-51.


\bibitem{FT1} R. Fuhrmann and F. Torres, {\em On Weierstrass points and
optimal curves}, Rend. Circ. Mat. Palermo {\bf 51}, (1998), 25--46.

\bibitem{GV} A. Garcia and J.F. Voloch, {\em Wronskians and independence
in fields of prime characteristic}, Manuscripta Math. {\bf 59}(1987),
457-469.


\bibitem{GV3} A. Garcia and J.F. Voloch, {\em Duality for projective
curves}, Bol. Soc. Brazil. Mat (N.S.) {\bf 21}, (1991) 159-175.

\bibitem{GeV} G. van der Geer and M. van der Vlugt, Generalized
Reed-M\"uller codes and curves with many points, preprint.
\bibitem{He} A. Hefez, {\em Non-reflexive curves}, Compositio
Math. {\bf 69} (1989), 3-35.

\bibitem{HK} A. Hefez and N. Kakuta, {\em Tangent envelopes of
higher order duals of projective curves}, Rend. Circ. Mat.
Palermo, Suppl. {\bf 51} (1998), 47-56.

\bibitem{HV} A. Hefez and J.F. Voloch, {\em Frobenius non
classical curves}, Arch. Math. {\bf 54} (1990), 263-273.



\bibitem{HT} J.W.P. Hirschfeld and J.A. Thas, ``General Galois geometries,''
Oxford University Press, Oxford, 1991.

\bibitem{HM} M. Homma, {\em Space curves with degenerate strict
duals}, Comm. Algebra {\bf 20} 867-874.

\bibitem{K} H. Kaji, {\em Strangeness of higher order space
curves}, Comm. Algebra {\bf 20} (1992), 1535-1548.


\bibitem{Lang} S. Lang, ``Abelian Varieties", Interscience, New
York, 1959.

\bibitem{RJ} J. Rathmann,{\em The uniform position principle for
curves in characteristic $p$}, Math. Ann. {\bf 276} (1987),
565-579.

\bibitem{RS} H.G. R\"uck and  H. Stichtenoth, A characterization
of Hermitian function fields over finite fields, {\em J. Reine Angew.
Math.} {\bf 457} (1994), 185--188.

\bibitem{SB} B. Segre, {\em Forme e geometrie Hermitiane, con
particolare riguardo al caso finito}, Ann. Mat. Pura Appl. {\bf
70} (1965) 1-201.


\bibitem{SV} K.O. St\"{o}hr and J.F. Voloch, {\em Weierstrass points and
curves over finite fields}, Proc. London Math. Soc. {\bf 52} (1986), 1-19.


\bibitem{Tate} J. Tate, {\em Endomorphisms of abelian varieties over
finite fields}, Inventiones Math. {\bf 2} (1966), 134-144.

\end{thebibliography}
\end{document}